\documentclass[12pt]{article}
\usepackage{amsmath,amsfonts,amssymb,amscd}

\newtheorem{thm}{Theorem}[section]

\newtheorem{lem}[thm]{Lemma}
\newtheorem{cor}[thm]{Corollary}

\newtheorem{prop}[thm]{Proposition}

\makeatletter \@addtoreset{figure}{section} \makeatother
\makeatletter
\long\def\@makecaption#1#2{%
   \vskip 10\p@
   \setbox\@tempboxa\hbox{{#1}\ \ #2}%
   \ifdim \wd\@tempboxa >\hsize
       {#1}\ \ #2\par
   \else
       \hbox to\hsize{\hfil\box\@tempboxa\hfil}%
   \fi}
\makeatother

\def\qed{\hfill \rule{4pt}{7pt}}
\def\pf{\noindent {\it Proof.} }
\title{\bf Long rainbow path in properly edge-colored complete
graphs\footnote{Supported by NSFC No.10901035 and No.11371205.}}

\author{
\small  He Chen$^1$ and Xueliang Li$^2$\\
[2mm]
\small $^1$Department of Mathematics,\\
\small  Southeast University, Nanjing 210096, China\\
\small chenhe@seu.edu.cn\\
\small $^2$Center for Combinatorics and LPMC \\
\small Nankai University, Tianjin 300071, China \\
\small lxl@nankai.edu.cn\\
}

\begin{document}
\date{}
\maketitle
\begin{abstract}

Let $G$ be an edge-colored graph. A rainbow (heterochromatic, or
multicolored) path of $G$ is such a path in which no two edges have
the same color. Let the color degree of a vertex $v$ be the number
of different colors that are used on the edges incident to $v$, and
denote it to be $d^c(v)$. It was shown that if $d^c(v)\geq k$ for
every vertex $v$ of $G$, then $G$ has a rainbow path of length at
least $\min\{\lceil\frac{2k+1}{3}\rceil,k-1\}$. In the present
paper, we consider the properly edge-colored complete graph $K_n$
only and improve the lower bound of the length of the longest
rainbow path by showing that if $n\geq 20$, there must have a
rainbow path of length no less than $\displaystyle
\frac{3}{4}n-\frac{1}{4}\sqrt{\frac{n}{2}-\frac{39}{11}}-\frac{11}{16}$.
\\
[2mm] {\bf Keywords:} properly edge-colored graph, complete graph,
rainbow ( heterochromatic, or multicolored) path.\\
[2mm] {\bf AMS Subject Classification (2010)}: 05C38, 05C15

\end{abstract}

\section{Introduction}

We use Bondy and Murty \cite{B-M} for terminology and notation not
defined here and consider simple graphs only.

Let $G=(V,E)$ be a graph. By an {\it edge-coloring} of $G$ we mean a
function $C: E\rightarrow \mathbb{N} $, the set of natural numbers.
If $G$ is assigned such a coloring, then we say that $G$ is an {\it
edge-colored graph}. Denote the edge-colored graph by $(G,C)$, and
call $C(e)$ the {\it color} of the edge $e\in E$. We say that
$C(uv)=\emptyset$ if $uv\notin E(G)$ for $u,v\in V(G)$. For a
subgraph $H$ of $G$, we denote $C(H)=\{C(e) \ | \ e\in E(H)\}$ and
$c(H)=|C(H)|$. For a vertex $v$ of $G$, the {\it color neighborhood}
$CN(v)$ of $v$ is defined as the set $\{C(e)\ | \ e \mbox{ is
incident with }v\}$, the {\it color degree} $d^c(v)=|CN(v)|$. A
subgraph of $G$ is called {\it rainbow (heterochromatic, or
multicolored)} if any two edges of it have different colors. If $u$
and $v$ are two vertices on a path $P$, $uPv$ denotes the segment of
$P$ from $u$ to $v$, whereas $vP^{-1}u$ denotes the same segment but
from $v$ to $u$.

There are many existing publications dealing with the existence of
paths and cycles with special properties in edge-colored graphs. The
heterochromatic Hamiltonian cycle or path problem was studied by
Hahn and Thomassen \cite{H-T}, R\"{o}dl and Winkler (see
\cite{F-R}), Frieze and Reed \cite{F-R}, and Albert, Frieze and Reed
\cite{A-F-R}. In \cite{A-J-Z}, Axenovich, Jiang and Tuza gave the
range of the maximum $k$ such that there exists a $k$-good coloring
of $E(K_n)$ that contains no properly colored copy of a path with
fixed number of edges, no heterochromatic copy of a path with fixed
number of edges, no properly colored copy of a cycle with fixed
number of edges and no heterochromatic copy of a cycle with fixed
number of edges, respectively. In \cite{E-T-1}, Erd\"{o}s and Tuza
studied the heterochromatic paths in infinite complete graph
$K_\omega$. In \cite{E-T-2}, Erd\"{o}s and Tuza studied the values
of $k$, such that every $k$-good coloring of $K_n$ contains a
heterochromatic copy of $F$ where $F$ is a given graph with $e$
edges ($e<n/k$). In \cite{M-S-T}, Manoussakis, Spyratos and Tuza
studied $(s,t)$-cycle in $2$-edge colored graphs, where
$(s,t)$-cycle is a cycle of length $s+t$ and $s$ consecutive edges
are in one color and the remaining $t$ edges are in the other color.
In \cite{M-S-T-V}, Manoussakis, Spyratos, Tuza and Voigt studied
conditions on the minimum number $k$ of colors, sufficient for the
existence of given types (such as families of internally pairwise
vertex-disjoint paths with common endpoints, hamiltonian paths and
hamiltonian cycles, cycles with a given lower bound of their length,
spanning trees, stars, and cliques ) of properly edge-colored
subgraphs in a $k$-edge colored complete graph. In \cite{C-M-M},
Chou, Manoussakis, Megalaki, Spyratos and Tuza showed that for a
2-edge-colored graph $G$ and three specified vertices $x, y$ and
$z$, to decide whether there exists a color-alternating path from
$x$ to $y$ passing through $z$ is NP-complete. Many results in these
papers were proved by using probabilistic methods.

In \cite{A-J-Z}, Axenovich, Jiang and Tuza considered the local
variation of anti-Ramsey problem, namely, they studied the maximum
$k$ such that there exists a $k$-good edge-coloring of $K_n$
containing no heterochromatic copy of a given graph $H$, and denote
it by $g(n,H)$. They showed that for a fixed integer $k\geq 2$,
$k-1\leq g(n, P_{k+1})\leq 2k-3$, i.e., if $K_n$ is edge-colored by
a $(2k-2)$-good coloring, then there must exist a heterochromatic
path $P_{k+1}$, and there exists an a $(k-1)$-good coloring of $K_n$
such that no heterochromatic path $P_{k+1}$ exists.

In \cite{B-L}, the authors considered long heterochromatic paths in
general graphs with a $k$-good coloring and showed that if $G$ is an
edge-colored graph with $d^c(v)\geq k$ (color degree condition) for
every vertex $v$ of $G$, then $G$ has a heterochromatic path of
length at least $\lceil\frac{k+1}{2}\rceil$. In \cite{C-L-1,C-L-2},
we got some better bound of the length of longest heterochromatic
paths in general graphs with a $k$-good coloring.

In \cite{C-L-3}, we showed that if $|CN(u)\cup CN(v)|\geq s$ (color
neighborhood union condition) for every pair of vertices $u$ and $v$
of $G$, then $G$ has a heterochromatic path of length at least
$\lceil\frac{s+1}{2}\rceil$, and gave examples to show that the
lower bound is best possible in some sense.

In \cite{G-M}, Gy\'{a}rf\'{a}s and Mhalla showed that in any
properly edge-colored complete graph $K_n$, there is a rainbow path
with no less than $(2n+1)/3$ vertices. In \cite{C-L-2} we got a
better result, showing that in any edge-colored graph $G$, if for
every vertex of $G$ there are at least $k$ colors appear on it, then
the longest rainbow path in $G$ is no shorter than
$\lceil\frac{2k}{3}\rceil+1$.

\begin{thm}\label{2/3 k}\cite{C-L-2}
Let $G$ be an edge-colored graph. If $d^c(v)\geq k$ for every vertex
$v\in V(G)$, then $G$ has a heterochromatic path of length at least
$\min\{\lceil\frac{2k}{3}\rceil+1,k-1\}$.
\end{thm}

 In this paper, we will improve the bound in
\cite{G-M}, and show that a longest rainbow path in a properly
edge-colored $K_n$ is not shorter than $\displaystyle
\left(\frac{3}{4}-o(1)\right)n$.

\section{Propositions of a longest rainbow path}

Suppose $G$ is a properly edge-colored $K_n$, $P=v_0v_1v_2\cdots
v_l$ is one of the longest rainbow paths in $G$, and
$C(v_{i-1}v_i)=C_i$ ($i=1,2,\cdots,l$).

Suppose $l<n-2$ and $u$ is an arbitrary vertex which does not belong
to the path $P$. Then we can easily get the following proposition.

\begin{prop}\label{P1}
$C(v_0u)\in C(P)$, $C(v_lu)\in C(P)$.
\end{prop}
\pf Otherwise, $uv_0Pv_l$ or $uv_lP^{-1}v_0$ is a rainbow path of
length $l+1$, a contradiction. \qed

\begin{prop}\label{P2}
If $C(uv_i)\notin C(P)$, then $C(uv_{i-1})\in C(P)$, $C(uv_{i+1})\in
C(P)$.
\end{prop}
\pf Otherwise, $v_0Pv_{i-1}uv_iPu_l$ or $v_0Pv_iuv_{i+1}Pu_l$ is a
rainbow path of length $l+1$, a contradiction. \qed

\begin{prop}\label{P3}
If $C(uv_i)\notin C(P)$, then $\{C(v_0v_{i+1}),
C(v_lv_{i-1})\}\subset C(P)\cup C(uv_i)$.
\end{prop}
\pf Otherwise, $uv_iP^{-1}v_0v_{i+1}Pv_l$ or
$uv_iPv_lv_{i-1}P^{-1}v_0$ is a rainbow path of length $l+1$, a
contradiction. \qed

\begin{prop}\label{P4}
If $C(uv_i)\notin C{P}$, then $C(v_0v_l)\in C(P)\setminus
\{C_{i-1},C_i\}$.
\end{prop}
\pf Otherwise, $uv_iPv_lv_0Pv_{i-1}$ or
$uv_iP^{-1}v_0v_lP^{-1}v_{i+1}$ is a rainbow path of length $l+1$, a
contradiction. \qed

\begin{prop}\label{P5}
If $C(v_0v_i)\notin C(P)$, then $C(v_lu)\in
C(P)\setminus{C(v_{i-1}v_i)}$; if $C(v_lv_i)\notin C(P)$, then
$C(v_0u)\in C(P)\setminus{C(v_iv_{i+1})}$.
\end{prop}
\pf Otherwise, $v_{i-1}P^{-1}v_0v_iPv_lu$ or
$v_{i+1}Pv_lv_iP^{-1}v_0u$ is a rainbow of length $l+1$, a
contradiction. \qed

\begin{prop}\label{P6}
If $C(v_0v_i)\notin C(P)$, then $C(v_{i-1}u)\in
C(P)\cup{C(v_0v_i)}$; if $C(v_lv_i)\notin C(P)$, then
$C(v_{i+1}u)\in C(P)\cup{C(v_lv_i)}$.
\end{prop}
\pf Otherwise, $uv_{i-1}P^{-1}v_0v_iPv_l$ or
$uv_{i+1}Pv_lv_iP^{-1}v_0$ is a rainbow of length $l+1$, a
contradiction. \qed

With these propositions, we can give new lower bound of a longest
rainbow path. And we will do that separately in the following two
situations: the biggest rainbow cycle is of length $l+1$, and the
biggest rainbow cycle is of length less than $l+1$.

\section{A longest rainbow path has the same number of vertices as
a biggest rainbow cycle}

If the longest rainbow path has the same number of vertices as the
biggest rainbow cycle, then the biggest rainbow cycle is of length
$l+1$, and there exists a rainbow path $P=v_0v_1\cdots v_l$ such
that $C(v_0v_l)\notin C(P)$.

Then, we can easily get the following conclusion from Proposition
\ref{P4}.

\begin{lem}\label{L1}
If $C(v_0v_l)\notin C(P)$, then for an arbitrary $u\in V(G\setminus
P)$, $C(u,P)\in C(P)\cup C(v_0v_l)$.
\end{lem}

By using this Lemma, we can get one of our main conclusions.

\begin{thm}\label{T1}
If $n\geq 20$ and $C(v_0v_l)\notin C(P)$, then $\displaystyle
l\geq\frac{3}{4}n-1$.
\end{thm}
\pf We will prove it by contradiction. Suppose a longest rainbow
path in $G$ is of length $l< \frac{3}{4}n-1$. Then $|V(G)\setminus
V(P)|=n-l-1> \frac{n}{4}\geq 5$.

We can conclude by Lemma \ref{L1} that for any vertex $u\in
V(G)\setminus V(P)$, $C(u,P)\subseteq C(P)\cup C(v_0v_l)$. On the
other hand, $|V(P)|=|C(P)\cup C(v_0v_l)|=l+1$ and $G$ is a properly
edge-colored $K_n$. Therefore, $C(u,P)=C(P)\cup C(v_0v_l)$,
$C(G\setminus P)\cap \left(C(P)\cup C(v_0v_l)\right)=\emptyset$.

Since $P$ is one of the longest rainbow paths, by Proposition
\ref{P1}, there exist $2\leq i_1<i_2<\cdots<i_{n-2-l}<l$, $1\leq
j_1<j_2<\cdots<j_{n-2-l}<l-1$, such that

$\begin{array}{ll} & |\{C(v_0v_{i_1}),C(v_0v_{i_2}),
\cdots,C(v_0v_{i_{n-2-l}})\}|\setminus(C(P)\cup
\{C(v_0v_l)\})\\
= & |\{C(v_lv_{j_1}),C(v_lv_{j_2}),
\cdots,C(v_lv_{j_{n-2-l}})\}|\setminus(C(P)\cup \{C(v_0v_l)\})\\
= & n-l-2
\end{array}$\\
Additionally, $C(uv_{i_k-1})\neq C(v_0v_l)$, $C(uv_{j_k+1})\neq
C(v_0v_l)$, $k=1,2,\cdots, n-l-2$.

Let $I=\{i-1| C(v_0v_i)\notin C(P)\cup C(v_0v_l), 2\leq i\leq l-1
\}$, $J=\{j+1|C(v_jv_l)\notin C(P)\cup C(v_0v_l), 1\leq j\leq
l-1\}$.

Now we distinguish the following two cases:

{\bf Case 1. $I\cap J\neq\emptyset$.}

This implies that there exists some $t$ in $I \cap J$, i.e.,
$$\{C(v_0v_{t+1}), C(v_lv_{t-1})\}\cap (C(P)\cup
C(v_0v_l))=\emptyset$$

{\bf Case 1.1. $C(v_0v_{t+1})\neq C(v_lv_{t-1})$.}

Since $n-l\geq 4$ and $C(u,P)=C(P)\cup C(v_0v_l)$, there are no less
than $3$ colors which is not in $C(P)\cup C(v_0v_l)$ such that they
belong to the color set $C(u,V(G)\setminus V(P))$. Therefore, there
exist $u_1,u_2 \in V(G)\setminus V(P)$ such that $C(u_1u_2)\notin
C(P)\setminus \{C(v_0v_l),C(v_0v_{t+1}),C(v_{t-1}v_l)\}$.

By Lemma \ref{L1}, there exists some vertex $v\in V(P)$ such that
$C(u_1v)=C(v_0v_l)$, denote it by $v_{i_0}$. We can conclude from
Proposition \ref{P6} that $i_0\neq t$. Since
$C'=v_0v_{t+1}Pv_lv_{t-1}P^{-1}v_0$ is a rainbow cycle of length $l$
in which the color $C(v_0v_l)$ does not appear on it. Therefore,
$u_2u_1v_{i_0}C$ contains a rainbow path of length $l+1$, a
contradiction.

{\bf Case 1.2. $C(v_0v_{t+1})=C(v_lv_{t-1})\notin C(P)\cup
C(v_0v_l)$.}

First, we can conclude that $C(v_{t-1}u)\neq C(v_0v_l)$ for any
vertex $u \in V(G)\setminus V(P)$. Otherwise, suppose there exists
some $u\in V(G)\setminus V(P)$ such that $C(v_{t-1}u)=C(v_0v_l)$.
Since $|V(G)\setminus V(P)|
>5$, there exists a vertex $u_1\in V(G)\setminus(V(P)\cup \{u\})$
such that $C(uu_1)\notin C(P)\cup \{C(v_0v_l), C(v_0v_{t+1})\}$.
Therefore, $u_1uv_{t-1}P^{-1}v_0v_{t+1}Pv_l$ is a rainbow path of
length $l+1$, a contradiction.

Then, we will show that $t-1\notin I\cup J$.

If $t-1\in I$, i.e., $C(v_0v_t)\notin C(P)\cup \{C(v_0v_l),
C(v_0v_{t+1})\}$, $C'=v_0Pv_{t-1}v_lP^{-1}v_tv_0$ is a rainbow cycle
of length $l+1$ without color $C(v_0v_l)$. On the other hand, by
\ref{L1} there exists a vertex $u\in V(G)\setminus V(P)$ and a
vertex $v_{i_0}\in V(P)$ such that $C(uv_{i_0})=C(v_0v_l)$. Then
$uv_{i_0}C$ contains a rainbow path of length $l+1$, a
contradiction.

If $t-1\in J$, i.e., $C(v_{t-2}v_l)\notin C(P)\cup \{C(v_0v_l),
C(v_{t-1}v_l)\}$, $C'=v_0Pv_{t-2}v_lP^{-1}v_{t+1}v_0$ is a rainbow
cycle of length $l$ without color $C(v_0v_l)$. Since $|V(G)\setminus
V(P)|>5$, for any vertex $u\in V(G)\setminus V(P)$, $d^c_{G\setminus
P}(u)\geq 5$. So, by Theorem \ref{2/3 k} there exists a rainbow path
$u_1u_2u_3\in G\setminus P$ with no colors in $C(P)\cup \{C(v_0v_l),
C(v_0v_{t+1}), C(v_{t-2}v_l)\}$. Since $G$ is properly edge-colored,
at least one edge in $\{v_tu_1, v_tu_3\}$ does not have color
$C(v_0v_l)$, W.O.L.G., assume $C(v_tu_1)\neq C(v_0v_l)$. Then,
because $C(v_{t-1}u_1)\neq C(v_0v_l)$, $C(u_1,P)=C(P)\cup
C(v_0v_l)$. So, by Lemma \ref{L1} there exists some $i_0$, $0\leq
i_0\leq l$, $i_0\neq t-1,t$ such that $C(u_1v_{i_0})=C(v_0v_l)$.
Then $u_3u_2u_1v_{i_0}C'$ contains a rainbow path of length $l+1$, a
contradiction.

So, we have $t-1\notin I\cup J$.

Let $K=I\cap J$, $I'=(I\setminus K)\cup \{t-1|t\in K\}$. Then
$|I'|=|I|$ and $I'\cap J =\emptyset$.

Additionally, for any $t\in I'\cup J$ and any $u\in V(G)\setminus
V(P)$, $C(v_tu)\neq C(v_0v_l)$. Otherwise, there exist some $t_0\in
K$ and some vertex $u\in V(G)\setminus V(P)$, such that
$C(v_{t_0-1}u)=C(v_0v_l)$. Since $|V(G)\setminus V(P)|\geq 6$, there
exists some vertex $u_1\in V(G)\setminus V(P)$ such that
$C(uu_1)\notin C(P)\cup \{C(v_0v_l), C(v_0v_{t_0+1})\}$. Then
$u_1uv_{t_0-1}P^{-1}v_0v_{t_0+1}Pv_l$ is a rainbow path of length
$l+1$, a contradiction.

On the other hand,
$$|I'\cup J|=|I'|+|J|=|I|+|J|\geq 2[(n-1)-(l+1)]=2(n-l-2),$$\\
and $|V(G)\setminus V(P)|=n-(l+1)=n-l-1$. So there are at least
$n-l-1$ $i$'s ($1\leq i\leq l-1$) such that $C(uv_i)=C(v_0v_l)$ for
some $u\in V(G)\setminus V(P)$. So we have $|I'\cup J|+(n-l-1)\leq
l-1$, and then $2(n-l-2)+n-l-1\leq l-1$, which implies $l\geq
\frac{3}{4}n-1$, a contradiction.

{\bf Case 2. $I\cap J=\emptyset$.}

By Proposition \ref{P6}, we have that for any $t\in I\cup J$ and any
$u\in V(G)\setminus V(P)$, $C(v_tu)\neq C(v_0v_l)$. On the other
hand, there are at least $|V(G)\setminus V(P)|=n-l-1$ $i$'s ($1\leq
i\leq l-1$) such that $C(uv_i)=C(v_0v_l)$ for some $u\in
V(G)\setminus V(P)$. So we have $|I\cup J|+(n-l-1)\leq l-1$, and
then $2(n-l-2)+n-l-1\leq l-1$, which implies $l\geq \frac{3}{4}n-1$,
a contradiction.

This complete the proof.\qed

\section{A biggest rainbow cycle has less vertices than a longest rainbow path}

Since a biggest rainbow cycle have less vertices than a longest
rainbow path, then $C(v_0v_l)\in C(P)$.

For any longest rainbow path $P$, by Proposition \ref{P1} and
Theorem \ref{T1}, there exist $2\leq i_1<i_2<\cdots<i_{t_1}<l$
($t_1\geq n-1-l$) such that
$$|\{C(v_0v_{i_1}), C(v_0v_{i_2}),\cdots, C(v_0v_{i_{t_1}})\}|= |CN(v_0)\setminus
C(P)|=t_1.$$

Now we will distinguish two cases: the case when there is a vertex
$u\in V(G)\setminus V(P)$ such that $C(v_lu)=C_1$, and the case when
there is no such vertex.

We first consider the case when there is a vertex $u\in
V(G)\setminus V(P)$ such that $C(v_lu)=C_1$.

\begin{thm}\label{T2}
If $C(v_0v_l)\in C(P)$ and there is a vertex $u\in V(G)$ such that
$C(v_lu)=C_1$, then $l\geq \displaystyle
\frac{3}{4}n-\frac{1}{4}\sqrt{\frac{n}{2}-\frac{39}{11}}-\frac{11}{16}$.
\end{thm}
\pf Suppose $P$ is a longest rainbow path that has the minimized
$t_1$.

We can conclude from Proposition \ref{P5} that $C_{i_k}\notin
C(v_l,C(G)\setminus V(P))$, $k=1,2,\cdots, t_1$.

Let $C^1=\{C_{i_k}| k=1,2,\cdots, t_1\}$,
$C_j^0=CN(v_{i_j-1})\setminus(C(P)\cup C(v_0v_{i_j}))$. Let the
color set $C_j^1$, $C_j^*$ ($j=1,2,\cdots,t_1$) be defined by the
following procedure.

For $j=1$ to $t_1$ do

\hspace*{10pt} $C_j^*=\emptyset$,

\hspace*{10pt} for $s=1$ to $i_j-3$

\hspace*{20pt} if $C(v_{i_j-1}v_s)\in C_j^0$, let
$C_j^*=C_j^*\cup\{C_{s+1}\}$;

\hspace*{10pt} for $s=i_{j+1}$ to $l-1$

\hspace*{20pt} if $C(v_{i_j-1}v_s)\in C_j^0$, let
$C_j^*=C_j^*\cup\{C_s\}$,

\hspace*{10pt} $C_j^1=C_{j-1}^1\cup C_j^*$

Then we can conclude that $|C_j^*|=|C_j^0|\geq t_1-1$ by Proposition
\ref{P1}.

Suppose $|C_{t_1}^1|-|C^0|=j_0$ and $j\geq j_0+2$.

Let $C_{j,1}=\{C(v_{i_t-1}v_{i_t})| t>j \mbox{~ and~}
C(v_{i_j-1}v_{i_t})\in C_j^0\}$,

\qquad $C_{j,2}=\{C(v_{i_t-1}v_{i_t})|t<j \mbox{~ and ~}
C(v_{i_j-1}v_{i_t-1})=C(v_0v_{i_t})\}$,

\qquad $C_{j,3}=\{C(v_{i_t-1})| t<j \mbox{~ and ~}
C(v_{i_j-1}v_{i_t-1})\in C_j^1\setminus C(v_0v_{i_t})\}$.

Then $C_{j,1}$, $C_{j,2}$, $C_{j,3}$ are mutually independent and
$C_j^*\cap C^0=C_{j,1}\cup C_{j,2}\cup C_{j,3}$. By the definition
$|C_{j,1}|\leq t_1-j$, $\displaystyle \bigcup\limits_{j=j_0+2}^{t_1}
C_{j,2}\subseteq \{C_{i_1}, C_{i_2},\cdots, C_{i_{t_1-1}}\}$ and
$C_{j,2}\cap C_{j',2}=\emptyset$ since $G$ is properly edge-colored.

Since $C(v_lu)=C_1$, we have $C_{j,3}=\emptyset$; otherwise,
$v_2Pv_{i_t-1}v_{i_j-1}P^{-1}v_{i_t}v_0v_{i_j}Pv_lu$ is a rainbow
path of length $l+1$, a contradiction. Therefore, $C_j^*\cap
C^0=C_{j,1}\cup C_{j,2}$.

On the other hand, $|C_j^*\setminus C^0|\leq |C_j^1\setminus
C^0|\leq j_0$. So, $|C_{j,2}|=|C_j^*\cap C^0|-|C_{j,1}|\geq
(t_1-1-j_0)-(t_1-j)=j-j_0-1$. Notice that $\displaystyle
\sum\limits_{j=j_0+2}^{t_1}|C_{j,2}|=\left|\bigcup\limits_{j=j_0+2}^{t_1}C_{j,2}
\right|\leq t_1-1$. Then, we have $\displaystyle
\sum\limits_{j=j_0+2}^{t_1} (j-j_0-1)\leq t_1-1$, i.e.,
$\displaystyle \frac{1}{2} \left( t_1^2 -2j_0t_1 -t_1 + j_0^2 +j_0
\right) \leq t_1-1$. Therefore, $\displaystyle j_0\geq
t_1-\frac{1}{2}-\sqrt{2t_1-\frac{7}{4}}$, $\displaystyle
|C_{t_1}|=t_1+j_0\geq 2t_1-\frac{1}{2}-\sqrt{2t_1-\frac{7}{4}}$.

Since $C(v_l, V(G)\setminus V(P))\subseteq C(P)\setminus
(C_{t_1}^1\cup \{C_l\})$ and $G$ is properly edge-colored,
$\displaystyle |V(G)\setminus V(P)|\leq
l-\left(2t_1-\frac{1}{2}-\sqrt{2t_1-\frac{7}{4}}\right) -1$, i.e.,
$\displaystyle n-(l+1)\leq
l-\left(2t_1-\frac{1}{2}-\sqrt{2t_1-\frac{7}{4}}\right) -1$. So,
$\displaystyle 2t_1-\sqrt{2t_1-\frac{7}{4}}\leq 2l-n+\frac{1}{2}$.
Since $f(x)=2x-\sqrt{2x-\frac{7}{4}}$ increases when $x>2$ and
$t_1\geq n-l-1>2$, we have
$$2(n-l-1)-\sqrt{2(n-l-1)-\frac{7}{4}}\leq 2l-n+\frac{1}{2}.$$
Therefore, $\displaystyle l\geq \frac{3}{4}n -
\frac{1}{4}\sqrt{\frac{n}{2}-\frac{39}{16}}-\frac{11}{16}$.

This completes the proof. \qed

Now we consider the case when for any longest rainbow path
$P=v_0v_1v_2\cdots v_l$ and any $u\in V(G)\setminus V(P)$,
$C(v_lu)\neq C_1$.

\begin{lem}\label{L2}
If for any longest rainbow path $P=v_0v_1v_2\cdots v_l$ and any
$u\in V(G)\setminus V(P)$, $C(v_lu)\neq C_1$ and there are at most
two $j$'s satisfying $2\leq j\leq t_1$, $i_j-i_{j-1}\geq 2$, then
$l\geq \frac{3n-4}{4}$.
\end{lem}
\pf For any $j$ ($1\leq j\leq t_1$), $v_{i_j-1}P^{-1}v_0v_{i_j}Pv_l$
is a rainbow path. So we can get by Proposition \ref{P5} and the
condition of this lemma that $\{C_{i_j-1}, C_{i_j}\}\cap C(v_l,
V(G)\setminus V(P))=\emptyset$.

Let $\displaystyle
C^*=\bigcup\limits_{j=1}^{t_1}\{C_{i_j-1},C_{i_j}\}$. Then
$|C^*|\geq 2t_1-2$ since there are at most two $j$'s satisfying
$2\leq j\leq t_1$, $i_j-i_{j-1}\geq 2$. On the other side, $C(v_l,
V(G)\setminus V(P))\subseteq C(P)\setminus (C^*\cup \{C_l\})$. So we
have
$$n-l-1\leq l-(2t_1-2)-1=l-2t_1+1\leq l-2(n-l-1)+1.$$
This implies that $\displaystyle l\geq \frac{3n-4}{4}$ and completes
the proof. \qed

Then we can get the following conclusion.

\begin{thm}\label{T3}
If $C(v_0v_l)\in C(P)$ and for any vertex $u\in V(G)$, $C(v_lu)\neq
C_1$, then $l\geq \displaystyle
\frac{3}{4}n-\frac{1}{4}\sqrt{\frac{n}{2}-\frac{39}{11}}-\frac{11}{16}$.
\end{thm}
\pf Let $i_0=\min\{i| \exists u\notin V(P) \mbox{~s.t.~}
C(v_lu)=C_i\}$. Suppose $P$ is one of the longest rainbow paths such
that $i_0$ is the smallest.

Let $j^*=\max\{j|i_j-i_{j-1}=1\}$. Then we have $i_0>i_{j^*}$;
otherwise, $v_1Pv_{i_{j^*-1}}v_0v_{i_{j^*}}Pv_l$ is also a rainbow
path of length $l$, but $C_{i_0}$ appears on the ($i_0-1$)-th edge
of the path, a contradiction.

Now we distinguish the following two cases.

{\bf Case 1. $i_0< i_{t_1}$.}

Let the integer $j_0$ and the color sets $C_j^0$, $C_j^*$,
$C_{j,1}$, $C_{j,2}$, $C_{j,3}$ be defined as in Theorem \ref{T2}.

Suppose $i_{j_1-1}<i_0<i_{j_1}$. Then we have that for any $j_1\leq
j_2\leq t_1$, $\{C(v_{j_2-1}v_{i_t-1})| 1\leq t<j_1\}\cap
C_{j_2}^0=\emptyset$. Otherwise, there exists $j_3<j_1\leq j_2$,
such that $C(v_{i_{j_3}-1}v_{i_{j_2}-1})\notin C_{j_2}^0$. Then,
$v_{i_{j_3}}Pv_{i_{j_2}-1}v_{i_{j_3}-1}P^{-1}v_0v_{i_{j_2}}Pv_l$ is
a rainbow path of length $l$, but the color $C_{i_0}$ appears on the
$(i_0-i_{j_3})$-th edge of this path, a contradiction to the choice
of $P$.

If there exists $j_1\leq j_2<j_3$ such that
$C(v_{i_{j_3}-1}v_{i_{j_2}-1})\notin \{C_{j_3}^0\cup
C(v_0v_{i_{j_2}})\}$, then
$v_1Pv_{i_{j_2}-1}v_{i_{j_3}-1}P^{-1}v_{i_{j_2}}v_0v_{i_{j_3}}Pv_l$
is a rainbow path of length $l$, but $C_{i_0}$ appears on the
$(i_0-1)$-th edge of this path, a contradiction.

Therefore, for any $j\geq j_1$, $C_{j,3}=\emptyset$,
$C_{j,2}\subseteq\{C_{i_t}|j_1\leq t<t_1\}$.

{\bf Case 1.1. $j_1>j_0$.}

As in Theorem \ref{T2}, we can get that $\displaystyle
\sum\limits_{j=j_1}^{t_1} (j-j_0-1)= \sum\limits_{j=j_1}^{t_1}
|C_{j,2}| = \left|\bigcup\limits_{j=j_1}^{t_1} C_{j,2}\right|\leq
t_1-j_1$. This implies that $(t_1-j_1+1)(j_1+t_1-2j_0-2)\leq
2(t_1-j_1)$. Therefore, $j_0\geq
\frac{1}{2}[(t_1^2-3t_1)-(j_1^2-3j_1)+2j_1-2]>\frac{1}{2}(2j_1-2)=j_1-1$,
a contradiction.

{\bf Case 1.2. $j_1\geq j_0$.}

By the same calculation we did in Theorem \ref{T2}, we can conclude
that $l\geq \displaystyle
\frac{3}{4}n-\frac{1}{4}\sqrt{\frac{n}{2}-\frac{39}{11}}-\frac{11}{16}$.

{\bf Case 2. $i_0>i_{t_1}$.}

If there are at most two $j$'s satisfying $2\leq j\leq t_1$,
$i_j-i_{j-1}\geq 2$, then by Lemma \ref{L2}, $\displaystyle l\geq
\frac{3n-4}{4}\geq \displaystyle
\frac{3}{4}n-\frac{1}{4}\sqrt{\frac{n}{2}-\frac{39}{11}}-\frac{11}{16}$.

So we will only consider the case when there are at least three
$j$'s satisfying $2\leq j\leq t_1$, $i_j-i_{j-1}\geq 2$. Then $l\geq
\frac{3n-4}{4}$. Suppose there are exactly $k$ ($k\geq 3$) such
$j$'s satisfying $s_1<s_2<\cdots<s_k$. Then for any integer $p$
($1\leq p\leq k$) $v_1Pv_{i_{s_p-1}}v_0v_{i_{s_p}}Pv_l$ is a rainbow
path of length $l$. Therefore, $$C(v_1,V(G)\setminus V(P))\subseteq
(C(P)\setminus\{C_{i_{s_p}}\})\cup \{C(v_0v_{i_{s_p-1}}),
C(v_0v_{i_{s_p}})\}.$$ Notice that $k\geq 3$, and so $\displaystyle
\bigcap\limits_{p=1}^k \{C(v_0v_{i_{s_p-1}}),
C(v_0v_{i_{s_p}})\}=\emptyset$, and then $C(v_1,V(G)\setminus
V(P))\subseteq (C(P)\setminus\{C_{i_{s_p}}\})$. Let
$C^*=C(P)\setminus \left( \bigcup\limits_{p=1}^s
C_{i_{s_p}}\cup\{C_1,C_2\}\right)$. Then $C(v_1,V(G)\setminus
V(P))\subset C^*$.

{\bf Case 2.1. $|C^*\cap \{C_1,C_2,\cdots, C_{i_{t_1}}\}|<t_1$.}

$|C^*\cap \{C_1,C_2,\cdots, C_{i_{t_1}}\}|<t_1$ implies that
$i_{t_1}-k-2<t_1$ and there exists a vertex $u\in V(G)\setminus
V(P)$ such that $C(v_1u)=C_t$, where $t\geq
i_{t_1}-[t_1-(i_{t_1}-k-2)]=t_1+k+2$ and it appears on the
$(l-t+1)$-th edge of the rainbow path
$v_lP^{-1}v_{i_{s_1}}v_0v_{i_{s_1-1}}P^{-1}v_1$ of length $l+1$. By
the choice of $P$, we can conclude that $l-t+1\geq i_0>i_{t_1}$,
i.e., $t\leq l-i_{t_1}+1$. Remember that $i_{t_1}\geq 2t_1-k$ and
$t_1\geq n-l-1$, and so we have $t_1+k+2\leq t\leq l-i_{t_1}+1\leq
l-2t_1+k+1$, i.e., $l\geq3t_1+1\geq 3n-3l-2$, and therefore
$\displaystyle l\geq\frac{3n-2}{4}\geq \displaystyle
\frac{3}{4}n-\frac{1}{4}\sqrt{\frac{n}{2}-\frac{39}{11}}-\frac{11}{16}$.

{\bf Case 2.2. $|C^*\cap \{C_1,C_2,\cdots, C_{i_{t_1}}\}|\geq t_1$.}

Suppose $C_t$ is the $t_1$-th color in $C^*$, $i_{j_0-1}< t\leq
i_{j_0}$ and there are $k_1$ $j$'s in the set $\{2,\cdots,j_0-1\}$
satisfying $i_j-i_{j-1}=1$. Then we can conclude that $t=t_1+k_1+2$
and $t> 2(j_0-1)-k_1-2=2j_0-k_1-4$. Since if $i_p-i_{p-1}>1$ then
$|C^*\cap\{C_{i_{p-1}+1},\cdots, C_{i_p}\}|\leq i_p-i_{p-1}$, we
have $$\displaystyle t_1=i_1-2+\sum\limits_{\begin{subarray}{c}
p\leq
j_0-1\\i_p-i_{p-1}>1\end{subarray}}(i_p-i_{p-1})+t-i_{j_0-1}=t_1-k_1-2.$$

On the other hand, $i_{t_1}\geq
i_{j_0}+2(t_1-j_0)-(k-k_1)=i_{j_0}+2t_1-2j_0-k+k_1\geq
t+2t_1-2j_0-k+k_1$. By Lemma \ref{L2} there is some integer $j$
satisfying $i_j-i_{j-1}=1$, and so
$v_lP^{-1}v_{i_j}v_0v_{i_{j-1}}P^{-1}v_1$ is a rainbow path of
length $l+1$ and $C_t$ appears on the $(l-t)$-th or $(l-t+1)$-the
edge. Therefore, we have $i_0\leq l-t$ by the choice of $P$. Then we
have $l-t\geq i_0>i_{t_1}\geq t+2t_1-2j_0-k+k_1$, i.e.,
\begin{eqnarray*}
l-t_1-k-2 &\geq& 3t_1+2k_1-2j_0+2\\ & >&
3t_1+2k_1+(-t-k_1-4)+2\\
&= &3t_1-t+k_1-2\\
&= &3t_1-(t_1+k_1+2)+k_1-2\\
&= &2t_1-4
\end{eqnarray*}
So, $l\geq 3t_1+k-2\geq 3t_1-2\geq 3(n-l-1)-2$, which implies that
$$\displaystyle l\geq \frac{3n-5}{4} \geq \displaystyle
\frac{3}{4}n-\frac{1}{4}\sqrt{\frac{n}{2}-\frac{39}{11}}-\frac{11}{16}.$$

This completes the proof. \qed

\section{Conclusion}

By Theorems \ref{T1}, \ref{T2} and \ref{T3}, we can easily get the
following conclusions.
\begin{thm}
For any properly edge-colored complete graph $K_n$ ($n\geq 20$),
there is a rainbow path of length no less than $\displaystyle
\frac{3}{4}n-\frac{1}{4}\sqrt{\frac{n}{2}-\frac{39}{11}}-\frac{11}{16}$.
\end{thm}

\begin{cor}
For any properly edge-colored complete graph $K_n$ ($n\geq 20$),
there is a rainbow path of length no less than $\displaystyle
(\frac{3}{4}-o(1))n$.
\end{cor}

\end{document}